\documentclass[12pt]{amsart} \headheight=8pt     \topmargin=0pt 
\def\cal{\mathcal} 

\usepackage{amsmath,amssymb,amsthm,latexsym}
\usepackage{amsfonts,array,amscd}

\input epsf

\def\epsfigbracket[#1=#2]{\epsffile{ps/#2}}
\def\epsfig#1{\epsfigbracket[#1]}
 
\theoremstyle{plain} 
\newtheorem{lemma}{Lemma}
\newtheorem{prop}[lemma]{Proposition} 
\newtheorem{theorem}[lemma]{Theorem} 
 
\newtheorem{fact}[lemma]{Fact} 
\newtheorem{keylemma}[lemma]{Keylemma} 

\theoremstyle{definition}  

\theoremstyle{remark}

\newcommand{\Z}{{\ensuremath{\mathbb Z}}}
\newcommand{\Q}{{\ensuremath{\mathbb Q}}}

\newcommand{\Aarhus}{\text{\sl \AA}}
\def\Zed{{\cal Z}}
\def\Ker{{\rm Ker}}
\def\lk{{\rm lk}}
\def\Tor{{\rm Tor}}
\def\cZ {{\check{Z}}^\sigma}
\newcommand\Zc{{\check{Z}}}
\def\exp{{\rm exp}}

\newcommand\strutu[2]{\,_{#1}\!\!\smile_{#2}}
\newcommand\strutn[2]{\,^{#1}\!\!\frown^{#2}}
\newcommand\A{{\cal A}}
\newcommand\wh{{\rm wh}}
\newcommand\stru{{\rm strut}}
\newcommand\defr{{\rm defr}}
\newcommand\Pwh{{P_{\rm wh}}}
\newcommand\Pstrut{{P_{\rm strut}}}
\newcommand\gl{{\rm gl}}
\newcommand\disuni{{\sqcup}}
\newcommand\Disuni{{\bigsqcup}}
\newcommand\source{{\rm source}}
\newcommand\target{{\rm target}}
\newcommand\chib{{\bar{\chi}}}

\newcommand\rank{{\rm rank}}
\newcommand\psib{{\bar{\psi}}}
\def\wheel{{
\setlength{\unitlength}{1.4192mm}
\begin{picture}(2,2.5)(-1,-1)
\thinlines \put(0,0){\circle{1.4142135}} 
\qbezier(0.5,0.5)(0.75,0.75)(1,1) 
\qbezier(0.5,-0.5)(0.75,-0.75)(1,-1)
\qbezier(-0.5,0.5)(-0.75,0.75)(-1,1) 
\qbezier(-0.5,-0.5)(-0.75,-0.75)(-1,-1) 
\end{picture}}}

\def\diagtheta{{ \setlength{\unitlength}{1.8192mm} 
\begin{picture}(2,2.5)(-1,-1) \thinlines 
\put(0,0){\circle{1.4142135}} \qbezier(0.5,0)(0,0)(-0.5,0) 
\end{picture}}} 
 
\def\isochord{{ \setlength{\unitlength}{1.8192mm} 
\begin{picture}(2,2.5)(-1,-1) \thicklines 
\put(0,0){\circle{1.4142135}} 
\thinlines\qbezier(0.5,0)(0,0)(-0.5,0) \end{picture}}}

\begin{document} { \parindent0cm 
\title[The LMO-invariant and the Alexander polynomial]{The 
LMO-invariant of $3$-manifolds of rank one and the Alexander 
polynomial} 

\date{} 

\author{Jens Lieberum} \address{\hskip-\parindent 
        Jens Lieberum\\
        Mathematical Sciences Research Institute\\
        1000 Centennial Drive\\
        Berkeley, CA 94720 - 5070\\
        USA}
\email{lieberum@msri.org}

\thanks{I would like to thank D.\ Bar-Natan, A.\ Beliakova, N.\ 
Habegger, T.\ Q.\ T.\ Le, and D.\ Thurston for helpful 
discussions. I thank the German Academic Exchange Service for 
financial support. 
Research at MSRI is supported in part by NSF grant DMS-9701755.} 


\begin{abstract} {\parindent0cm We prove that the LMO-invariant of 
a $3$-manifold of rank one is determined by the Alexander 
polynomial of the manifold, and conversely, that the Alexander 
polynomial is determined by the LMO-invariant. Furthermore, we 
show that the Alexander polynomial of a null-homologous knot in a 
rational homology $3$-sphere can be obtained by composing the 
weight system of the Alexander polynomial with the \AA rhus 
invariant of knots. 

\bigskip

 {Mathematics Subject Classification (2000): 57M25, 
57N65, 57M15. 

} 

} \end{abstract} 


\maketitle 

\section*{Introduction}

In analogy with the theory of Vassiliev invariants of links, 
different notions of finite type invariants of $3$-manifolds have 
been introduced. For integral homology spheres these different 
notions coincide with the original definition of Ohtsuki 
(\cite{Oht}). The LMO-invariant~$Z^{LMO}$ assembles 
all~$\Q$-valued finite type invariants of integral homology 
spheres in a formal series and is therefore called a universal 
finite type invariant (\cite{LMO}, \cite{Le}). For connected 
closed manifolds~$M$ the following is known about~$Z^{LMO}$: 

$$
\begin{array}{|l|l|}

\hline 

\mbox{For $M$ with ...} & \mbox{$Z^{LMO}$ is determined by and 
determines ... }\\ 

\hline \hline

 \rank\; H_1(M)=0 & \mbox{all $\Q$-valued invariants of Goussarov and Habiro (\cite{Habi})} \\ 

\hline 

H_1(M)=\Z & \mbox{the Alexander polynomial (Theorem 1 
of~\cite{GaH})}\\ 

\hline

\rank\; H_1(M)\geq 2 & \mbox{the Casson-Walker-Lescop invariant 
(\cite{HaB}, \cite{Hab}, \cite{Les})}\\ 

\hline 
\end{array}
$$

In this article we fill in the missing puzzle piece for the 
interpretation of the LMO-invariant of manifolds of rank~$\geq 1$ 
in terms of classical invariants. We prove the following 
generalization of Theorem~1 of~\cite{GaH}. 

\begin{theorem}\label{t:ZlmoNabla}
Let~$M$ be a closed oriented~$3$-manifold of rank~$1$. Then the 
LMO-invariant~$Z^{LMO}(M)$ is determined by the Alexander 
polynomial~$\nabla(M)$, and conversely,~$\nabla(M)$ is determined 
by~$Z^{LMO}(M)$. 
\end{theorem}

In the proof of Theorem~1 of~\cite{GaH} it was used that the 
Alexander polynomial~$\nabla$ of links~$L$ in~$S^3$ can be 
obtained from the universal Vassiliev invariant~$Z$ of links 
in~$S^3$ via a map~$W_\nabla$ as follows:

\begin{equation}\label{e:nablacan}
\frac{h}{e^{h/2}-e^{-h/2}}\nabla(L)_{\vert\;t^{1/2}=e^{h/2}}= 
W_\nabla\circ Z(L). 
\end{equation}

We generalize Equation~(\ref{e:nablacan}) by replacing~$Z$ by the 
\AA rhus invariant~$\Aarhus$ (\cite{BGRT1}) of knots in a rational 
homology sphere: 

\begin{theorem}\label{t:nablacan} Let~$K$ be a null-homologous 
knot in a rational homology $3$-sphere. Then 
$$\frac{h}{e^{h/2}-e^{-h/2}}\nabla(K)_{\vert\;t^{1/2}=e^{h/2}}= 
W_\nabla\circ\Aarhus(K).$$ \end{theorem} 

Theorem~\ref{t:nablacan} is an important ingredient in the proof 
of Theorem~\ref{t:ZlmoNabla}. Theorems~\ref{t:ZlmoNabla} 
and~\ref{t:nablacan} will be proven in Section~\ref{s:proofs}. In 
Sections~\ref{s:alex}--\ref{s:univ} we prepare these proofs by 
recalling definitions and properties of the Alexander polynomial 
of links and manifolds, of unitrivalent diagrams and the 
map~$W_\nabla$, and of the universal finite type 
invariants~$Z$,~$Z^{LMO}$ and~$\Aarhus$. 
 

 \section{The Alexander polynomial}\label{s:alex}

In this section we make preliminary definitions and recall some 
facts about the Alexander polynomial~$\nabla$ from~\cite{Les}. 

All manifolds and submanifolds in this paper are oriented. Let~$M$ 
be a rational homology $3$-sphere (meaning that~$M$ is a connected 
closed manifold of dimension~$3$ with~$H_1(M,\Q)=0$). Let~$K$ be a 
knot in~$M$. Choose a tubular neighborhood~$T$ of~$K$. A meridian 
of~$K$ is a simple closed curve~$m$ on the boundary~$\partial T$ 
of~$T$ that is null-homologous in~$T$. The curve~$m$ is oriented 
by the right-hand rule. There exists a unique isomorphism 
$i_K:H_1(M\setminus K,\Q)\longrightarrow \Q$ that sends a meridian 
of~$K$ to~$1$. As a $\Q$-linear map~$i_K$ is uniquely determined 
by the property that for any oriented surface~$\Sigma\subset M$ 
with~$\partial\Sigma\cap K=\emptyset$ the 
value~$i_K(\partial\Sigma)$ is the intersection number of~$K$ 
with~$\Sigma$. 
 For disjoint knots $K_1, K_2$ the linking number~$\lk(K_1, K_2)$ 
is defined as~$i_{K_1}(K_2)$. The linking number~$\lk(.,.)$ is 
symmetric. 

Denote the number of components of a link~$L$ by~$\vert L\vert$. A 
framed link~$L$ is a link with a simple closed curve~$\mu_i$ on 
the boundary~$\partial T_i$ of a tubular neighborhood~$T_i$ of 
each component~$K_i$ ($i=1,\ldots,\vert L\vert$). Inside of~$T_i$, 
$\mu_i$ is homologous to~$q_i K_i$ for some~$q_i\in\Z$. The {\em 
linking matrix}~$(l_{ij})$ of~$L$ is defined 
by~$l_{ij}=\lk(K_i,\mu_j)/q_j$. The link~$L$ has integral framing 
if all~$q_i$ are~$1$. The values~$l_{ii}\in\Q$ are called framing 
of~$K_i$. We denote by~$M_L$ the manifold obtained by surgery 
on~$L\subset M$. 

Let~$L\subset M$ be a null-homologous link. Then there exists an 
oriented connected surface~$\Sigma$ embedded in~$M$ such 
that~$\partial \Sigma=L$. Any surface with this property is called 
a Seifert surface of~$L$. Let $\Sigma^\pm=\Sigma^+\cup\Sigma^-$ be 
a tubular neighborhood of~$\Sigma$ such 
that~$\Sigma=\Sigma^+\cap\Sigma^-$ and~$\Sigma^+$ lies on the 
positive side of~$\Sigma$. The Seifert form of~$\Sigma\subset M$ 
is the $\Z$-bilinear form $s:H_1(\Sigma)\times 
H_1(\Sigma)\longrightarrow\Q$ defined by sending homology 
classes~$a$, $b$ to $\lk(A^-,B^+)$ where~$A^-$ is a knot 
in~$\Sigma^-$ representing~$a$ and~$B^+$ is a knot in~$\Sigma^+$ 
representing~$b$. In this section a matrix of~$s$ with respect to 
an arbitrary basis of~$H_1(\Sigma)$ is called a Seifert matrix 
of~$\Sigma$ (later we will choose a particular basis 
of~$H_1(\Sigma)$). Define the bilinear form $s^*$ by 
$s^*(a,b)=s(b,a)$. Then $s-s^*$ is the intersection form 
of~$\Sigma$. Denote the transpose of a matrix~$V$ by~$V^*$. 

\begin{prop}\label{p:nablaKS} Let~$L$ be a null-homologous link in 
a rational homology sphere~$M$. Choose a Seifert surface~$\Sigma$ 
of~$L$. Let~$V$ be a Seifert matrix of~$\Sigma$. Then 
$$\nabla(L)=\det 
(t^{1/2}V-t^{-1/2}V^*)\in(t^{1/2}-t^{-1/2})^{\vert L\vert-1} 
\Q[(t^{1/2}-t^{-1/2})^2]\subset\Q[t^{\pm 1/2}] $$ 
 
is an invariant of the pair $L\subset M$ up to homeomorphism; in 
particular it is an isotopy invariant of~$L$. 
\end{prop}

Proposition~\ref{p:nablaKS} can be proven by using sign-determined 
Reidemeister torsion (see Proposition~2.3.13 of~\cite{Les}, 
\cite{Tur}). 

Up to sign the invariant~$\nabla(L)$ can be described as follows. 
Let~$N$ be a connected $3$-manifold and let 
$\varphi:H_1(N)\longrightarrow\Zed=\Z$ be a homomorphism. Let 
$\widetilde{N}$ be the connected cover of~$N$ corresponding to 
$\Ker(\varphi)$. Then~$H_1(\widetilde{N})$ is a module over the 
group ring~$\Z[\Zed]\cong\Z[t^{\pm 1}]$. Let~$J\subset\Z[t^{\pm 
1}]$ be the order ideal of~$H_1(\widetilde{N})$. Let 
$\Delta_\varphi(N)$ be a generator of the smallest principal ideal 
containing~$J$. Then~$\Delta_\varphi(N)$ is unique up to 
multiplication by~$\pm {t^i}$. For a link in a rational homology 
sphere~$M$ we denote~$\Delta_\varphi(M\setminus L)$ 
by~$\Delta(L)$, where $\varphi:H_1(M\setminus L)\longrightarrow 
\Z$ is given by the sum of the linking numbers with the components 
of~$L$. The following lemma (see Proposition~2.3.13 of~\cite{Les}) 
relates~$\nabla(L)$ and~$\Delta(L)$. 

\begin{lemma}\label{l:nabladelta} Let~$L$ be a null-homologous 
link in a rational homology sphere~$M$. Then there exists a 
unique~$i\in\Z$ such that $t^{i/2}\Delta(L)$ is invariant under 
the replacement of~$t^{1/2}$ by~$-t^{-1/2}$. For 
some~$\epsilon\in\{\pm 1\}$ we have $\epsilon 
t^{i/2}\Delta(L)=\vert H_1(M)\vert \nabla(L)$. \end{lemma} 

Now consider a connected closed $3$-manifold~$N$ of rank~$1$. 
Denote the quotient of~$H_1(N)$ by its torsion 
subgroup~$\Tor(H_1(N))$ by~$H_1^\#(N)$. Choose an 
isomorphism~$\psi:H_1^\#(N)\longrightarrow\Z$. Denote the 
composition of the canonical projection~$H_1(N)\longrightarrow 
H_1^\#(N)$ with~$\psi$ by~$\psib$. The following two lemmas 
(see~\cite{Les}, Section~5.1) allow to compare~$\Delta_{\psib}(N)$ 
with a knot invariant. 
 
\begin{lemma}\label{l:surga} Every connected closed 
$3$-manifold~$N$ of rank~1 can be obtained by $0$-framed surgery 
on a null-homologous knot~$K$ in a rational homology sphere~$M$. 
We then have $\Tor(H_1(N))\cong H_1(M)$. \end{lemma} 

\begin{lemma}\label{l:nablaKMK} Let~$K$ be a null-homologous 
$0$-framed knot in a rational homology $3$-sphere~$M$. 
Then~$\Delta(K)$ is equal to~$\Delta_{\psib}(M_K)$ up to 
multiplication by~$\pm t^i$. \end{lemma} 

We see by Lemmas~\ref{l:nabladelta}, \ref{l:surga} 
and~\ref{l:nablaKMK} that there exists $j\in\Z$ such that 
$t^{j/2}\Delta_{\psib}(N)$ is invariant under the replacement of 
$t^{1/2}$ by $-t^{-1/2}$. Furthermore, we can 
choose~$\epsilon\in\{\pm 1\}$ such that 
$\epsilon\Delta_{\psib}(N)_{\vert t=1}=\vert H_1(\Tor(N))\vert>0$. 
Denote~$(\epsilon t^{j/2}/\vert H_1(\Tor(N))\vert 
)\Delta_\psib(N)$ by~$\nabla(N)$. The definition of~$\nabla(N)$ 
does not depend on the choice of the isomorphism~$\psi$ 
because~$\nabla(N)\in\Q[(t^{1/2}-t^{-1/2})^2]$. The 
invariant~$\nabla$ satisfies 

\begin{equation}\label{e:nablaKMK}
\nabla(K)=\nabla(M_K). \end{equation}

for all null-homologous $0$-framed knots~$K$ in a rational 
homology sphere~$M$.

\section{Unitrivalent diagrams and $W_\nabla$}\label{s:Wnabla}

In this section we briefly recall facts about unitrivalent 
diagrams and use them to state properties of the Vassiliev 
invariants in the Alexander polynomial~$\nabla$. 

 Let~$\Gamma$ be a 
compact oriented $1$-manifold whose boundary~$\partial \Gamma$ is 
partitioned into two ordered sets called upper and lower boundary. 
Let~$X$ be a set. A {\em unitrivalent diagram} with {\em 
skeleton}~$\Gamma$ is a graph~$D$ with distinguished 
subgraph~$\Gamma$ such that all vertices of~$D$ are either 
univalent or trivalent. Trivalent vertices not lying on~$\Gamma$ 
are called {\em internal} and are oriented by a cyclic order of 
the incident edges. Univalent vertices are also called {\em legs}. 
Each leg of a unitrivalent diagram is labeled by an element 
of~$X$. We allow connected components in~$D$ that do not 
intersect~$\Gamma$ whenever these components contain at least one 
trivalent vertex. Recall the definition of a $\Q$-vector 
space~$\A(\Gamma,X)$ generated by unitrivalent diagrams modulo 
relations called~$(STU)$,~$(IHX)$, and~$(AS)$ (\cite{BN1}). 
When~$\Gamma$ is equipped with additional information (for 
example: dots on circle-components of~$\Gamma$, a set~$Y$ in 
bijection with circle-components of~$\Gamma$, a distinguished 
subset of the components of~$\Gamma$,...), we require in the 
definition of~$\A(\Gamma,X)$ that homeomorphisms between 
unitrivalent diagrams also preserve this additional data. The 
space~$\A(\Gamma, X)$ is graded by half of the number of vertices 
of unitrivalent diagrams. Denote~$\A(\Gamma,\emptyset)$ 
by~$\A(\Gamma)$. 

The invariants of $\ell$-component links in~$S^3$ that are 
coefficients of~$z^i=(t^{1/2}-t^{-1/2})^i$ in~$\nabla(L)$ induce 
linear forms $W_i:\A^\ell_i\longrightarrow\Q$ on the degree-$i$ 
part~$\A^\ell_i$ of~$\A^\ell:=\A({S^1}^{\disuni\ell})$ (see 
Section~3 of~\cite{BNG}) . For~$a$ in the completion of~$\A^\ell$ 
by the degree, we define 

$$W_\nabla(a)=\sum_i W_i(a)h^i\in\Q[[h]].$$

 It will follow from 
Theorem~\ref{t:nablacan} and can also be seen directly that the 
Alexander polynomial of links in a rational homology sphere 
induces the same map~$W_\nabla$ (see skein relation~2.3.16 
of~\cite{Les}, or Exercise~3.10 of~\cite{BNG}). The map~$W_\nabla$ 
and its extensions to~$\A(\Gamma,X)$ obtained from representations 
of the Lie superalgebra~$\gl(1\vert 1)$ have the following 
property (see Proposition~7.1 of~\cite{Vai}, consider the 
element~$\diagtheta$ of~$\A(\emptyset)$ seperately). 

\begin{lemma}\label{l:Wnabla}
Let $D\in\A(\Gamma,X)$ be a unitrivalent diagram. Assume that~$D$ 
has an internal vertex~$u$ such that all edges incident to~$u$ are 
connected to internal vertices. Then we have~$W_\nabla(D)=0$. 
\end{lemma}

Let $I_x\cong I:=[0,1]$ ($x\in X$). Denote the disjoint union 
by~$\disuni$. For every partition 
of~$\partial(\Gamma\disuni\Disuni_{x\in X}I_x)$ into two ordered 
sets called upper and lower boundary there exists an isomorphism 

\begin{equation}\label{e:defchi} 
\chi_X:\A(\Gamma,X\disuni 
Y)\longrightarrow\A\left(\Gamma\disuni\Disuni_{x\in X} 
I_x,Y\right)  
\end{equation}

given by the average over all permutations of putting~$x$-labeled 
univalent vertices of a diagram on the corresponding skeleton 
component~$I_x\cong I$ of~$\Gamma\disuni\Disuni_{x\in X} I_x$. The 
inverse of~$\chi_X$ will be denoted by~$\sigma_X$ and the set~$X$ 
will not be specified when it is clear from the context. 
Obviously, there exists an isomorphism of 
$\A(\Gamma\disuni\Disuni_{x\in X} I_x,Y)$ with a space 
$\A(\Gamma\disuni\Disuni_{x\in X} {S^1_x}^*,Y)$, where the 
circles~${S^1_x}^*$ have a dot and are in bijection with the 
set~$X$. Similarly, we have a surjective map from 
$\A(\Disuni_{x\in X} I_x\disuni\Disuni_{y\in Y} I_y,Z)$ to  
$\A(\Disuni_{x\in X} I_x\disuni\Disuni_{y\in Y} {S^1_y},Z)$ given 
by closing the intervals~$I_y$ to form the circles~$S^1_y$. Denote 
the composition of~$\chi_Y$ with this surjective map 
by~$\bar{\chi}_Y$. 

An important special case is 
$\A(S^1)\cong\A({S^1}^*,\emptyset)\cong\A(I)=:\A$ 
(see~\cite{BN1}). The space~$\A$ is a commutative algebra with 
multiplication~$\#$ induced by the connected sum of the 
skeletons~$S^1$ of diagrams (resp.\ by the concatenation of 
skeletons~$I$ of diagrams). More generally, the connected sum 
of~$S^1$ with any distinguished skeleton component~$C$ of a 
unitrivalent diagram turns~$\A(\Gamma\cup C, X)$ into 
an~$\A$-module. Let~$\bar{\A}$ be the quotient of~$\A$ by the 
ideal generated by the element~$\isochord$ and 
let~$\pi:\A\longrightarrow\bar{\A}$ be the canonical projection. 
There exists a unique inclusion of 
algebras~$i:\bar{\A}\longrightarrow\A$ with the property 
that~$i(D)=D$ for all diagrams~$D$ such that~$D\setminus S^1$ is 
connected and~$D$ contains an internal vertex (\cite{BN1}, 
Equation~(5), Exercise~3.16). The map 
$P_{\defr}=i\circ\pi:\A\longrightarrow\A$ is called {\em deframing 
projection}. 

The disjoint union of unitrivalent diagrams turns~$\A(\emptyset, 
X)$ into a commutative algebra and~$\A(\Gamma, X)$ into an 
$\A(\emptyset, X)$-module. Important examples of diagrams 
in~$\A(\emptyset,X)$ are so-called {\em struts}~$\strutn{i}{j}$ 
with labels $i,j\in X$, and so-called {\em 
wheels}~$\omega_n=\wheel$ having $n$~internal vertices lying on a 
circle and~$n$ univalent vertices with the same label ($n=4$ in 
this example). Let $\A(\emptyset,X)_\stru\subset\A(\emptyset,X)$ 
be the subalgebra generated by struts and~$\A(\emptyset,X)_\wh$ be 
the subalgebra generated by wheels. It is known that 
$\A(\emptyset,X)_\stru$ is a polynomial algebra in the $n (n+1)/2$ 
different struts ($n=\vert X\vert$) and~$\A(\emptyset,X)_\wh$ is a 
polynomial algebra in wheels with an even number of univalent 
vertices. There exist unique projections from~$\A(\Gamma,X)$ 
to~$\A(\emptyset,X)_\stru$ (resp.\ $\A(\emptyset,X)_\wh$) that 
send all diagrams to~$0$ that have a connected component that is 
not a strut (resp.\ a wheel). 
Define~$\Pstrut:\A(\Gamma\disuni\Disuni_{x\in X} 
I_x,\emptyset)\longrightarrow \A(\Gamma,X)$ as the composition 
of~$\sigma$ with the projection 
to~$\A(\emptyset,X)_\stru\subset\A(\Gamma,X)$.  The map~$\Pstrut$ 
descends to~$\A(\Gamma\disuni\Disuni_{x\in X} S^1_x,\emptyset)$ 
where the circle-components~$S^1_x$ are in bijection with~$X$. 
Define $\Pwh:\A\longrightarrow\A(\emptyset,\{x\})$ as the 
composition of~$\sigma\circ P_{\defr}$ with the projection 
to~$\A(\emptyset,\{x\})_\wh$. We have~$\Pwh(a\# 
b)=\Pwh(a)\disuni\Pwh(b)$ for all $a,b\in\A$. The map~$\Pwh$ is 
related to~$W_\nabla$ as follows (see~\cite{Vai},~\cite{Kri}). 
 
\nopagebreak{ \begin{lemma}\label{l:WnablaPwh} For~$D\in\A$ the 
value $W_\nabla(D)$ depends only on~$\Pwh(D)$ and is determined by  
 $$W_\nabla(D_1\# 
D_2)=W_\nabla(D_1)W_\nabla(D_2)\quad\mbox{and}\quad 
W_\nabla(\bar{\chi}(\omega_{2n}))=-2h^{2n}.$$ \end{lemma} 

Lemma~\ref{l:WnablaPwh} was used in proofs of the 
Melvin-Morton-Rozansky conjecture~(\cite{BNG}). }
 
\section{Universal finite type invariants}\label{s:univ}

Recall from Section~3 of~\cite{LM2} that a non-associative framed 
tangle (or q-tangle)~$T$ is a usual tangle with integral framing, 
except that~$\source(T)$ and~$\target(T)$ are equipped with 
parentheses on the sequences of $\pm$-symbols associated with the 
lower and upper boundary points of~$T$. We denote by~$Z$ the 
universal Vassiliev invariant of non-associative framed tangles 
(see~\cite{LM2}). Denote the underlying $1$-manifold of a 
tangle~$T$ (together with the partition of~$\partial T$ into two 
ordered sets and possibly together with a decoration of~$T$ such 
as dots, distinguished components, ...) by~$\Gamma(T)$. Then the 
values~$Z(T)$ lie in the completion of~$\A(\Gamma(T))$ by the 
degree. 

Let $\nu=Z(O)$ be the invariant of the trivial knot with 
$0$-framing. Let~$T=L'\cup T''$ be a diagram of a framed 
non-associative tangle where the components of the sublink~$L'$ 
of~$T$ are in bijection with a set~$X'$ and each component of~$L'$ 
has a dot on its circle. Define $\Zc(T)$ as the connected sum 
of~$Z(T)$ with $\nu^{\otimes \vert L'\vert}$ along the components 
of~$\Gamma(L')$. Cut the chord diagrams in~$\Zc(T)$ at the dots 
and apply the isomorphism~$\sigma_{X'}$. The result lies in the 
completion of~$\A(\Gamma(T''),X')$ and is called~$\cZ(T)$. The 
value~$\cZ(T)$ is not invariant under isotopies of the tangle 
represented by the diagram~$T$.  For tangles~$T$ with dotted 
circles~$L'$ invariants~$Z^{LMO}_0(T)$ and~$\Aarhus_0(T)$ of 
isotopy (that are also invariant under second Kirby moves 
along~$L'$) are obtained from~$C=\cZ(T)$ as follows 
(see~\cite{LMO}, \cite{Le2}, \cite{BGRT2}). 

{\em Definition of~$Z^{LMO}_0$}: The degree-$n$ part of 
$Z^{LMO}_0(T):=<C>$ is obtained from the degree $n+\vert L'\vert 
n$ part of~$C$ by forgetting the diagrams in~$C$ that do not have 
exactly~$2n$ legs of each color~$x\in X'$, by summing over all the 
$((2n-1)!!)^{\vert L'\vert}=((2n)!/2^nn!)^{\vert L'\vert}$ 
possible ways of gluing pairs of legs of diagrams in~$C$ with the 
same label and by replacing circles that do not belong 
to~$\Gamma(T'')$ by~$-2n$. 

{\em Definition of $\Aarhus_0$}: $\Aarhus_0$ is only defined when 
the linking matrix~$(l_{ij})$ of~$L'$ is invertible (or 
equivalently, when~$S^3_{L'}$ is a rational homology sphere). 
Write~$C$ in the form 

$$C=P\disuni \exp\left(\frac{1}{2}\sum_{i,j\in X'} l_{ij} 
\strutn{i}{j}\right)$$ 

where~$P$ contains no struts. Let $(l^{ij})$ be the inverse matrix 
of $(l_{ij})$. Then 

$$ \Aarhus_0(T):=<P,\exp \left(-\frac{1}{2}\sum_{ i,j\in X'} 
l^{ij} \strutu{\partial i}{\partial j}\right)>, $$ 

where $<D_1,D_2>$ is $0$ if for some~$i$ the number of $i$-labeled 
legs of~$D_1$ is not equal to the number of~$\partial i$-labeled 
legs of~$D_2$, and is given by the sum of all ways of gluing all 
legs with $i$-labels to legs with $\partial i$-labels in the 
remaining case. 

Let $L\subset M$ be a link in a $3$-manifold. Represent~$L\subset 
M$ by a diagram~$L'\cup L''$ of a link in~$S^3$, such that 
$S^3_{L'}\cong M$ and the image of~$L''$ in $S^3_{L'}$ is mapped 
to~$L$ by this homeomorphism. Put a dot on each component of~$L'$. 
Two invariants~$Z^{LMO}$ and~$\Aarhus$ of homeomorphisms of the 
pair $(M,L)$ are obtained from~$Z^{LMO}_0(L'\cup L'')$ and 
$\Aarhus_0(L'\cup L'')$ by normalization (making it invariant 
under the first Kirby move) as follows: 

\begin{eqnarray}
Z^{LMO}(L) & = & 
Z^{LMO}_0(U_+)^{-\sigma_+}Z^{LMO}_0(U_-)^{-\sigma_-}Z^{LMO}_0(L'\cup 
L''),\\ \Aarhus(L) & = & 
\Aarhus_0(U_+)^{-\sigma_+}\Aarhus_0(U_-)^{-\sigma_-}\Aarhus_0(L'\cup 
L''),
\end{eqnarray}

where $U_\pm$ is the trivial knot with a dot and framing~$\pm 1$ 
and~$\sigma_+$ (resp.~$\sigma_-$) is the number of positive  
(resp.\ negative) eigenvalues of the linking matrix~$(l_{ij})$ 
of~$L'$. The invariants of the empty link~$Z^{LMO}(\emptyset)$ 
and~$\Aarhus(\emptyset)$ are also denoted by~$Z^{LMO}(M)$ 
and~$\Aarhus(M)$, respectively. The series~$\Aarhus_0(U_\pm)$ have 
degree-$0$ term~$1$. Therefore Lemma~\ref{l:Wnabla} implies 

\begin{equation}\label{e:A0A}
 W_\nabla\circ\Aarhus_0(L'\cup L'')=W_\nabla\circ\Aarhus(L). 
\end{equation} 

We will make use of the following result of~\cite{BGRT3} 
(Equation~(\ref{e:ALMO}) follows from Proposition~1.2 
of~\cite{BGRT3} in the same way as Theorem~1 of~\cite{BGRT3}): 

\begin{equation}\label{e:ALMO}
\Aarhus(L)=\vert H_1(M)\vert^{-deg} Z^{LMO}(L), 
\end{equation}

where $\vert H_1(M)\vert^{-deg}$ denotes the operation of 
multiplying diagrams of degree~$m$ by $\vert H_1(M)\vert^{-m}$. 

Let us recall some notation used in Lemma~\ref{l:A0gen} below. 
Let~$T$ be a non-associative framed tangle~$T$ with a 
distinguished subset~$\tilde{T}$ of its components. Denote 
by~$d(T)$ the non-associative framed tangle given by replacing 
each component in~$\tilde{T}$ by two copies that are parallel with 
respect to the framing. The symbols~$a\in\{+,-\}$ in~$\source(T)$ 
(resp.~$\target(T)$) that belong to~$\tilde{T}$ are replace by 
$(a\; a)$ in~$\source(d(T))$ (resp.~$\target(d(T))$). 
Define~$s(T)$ by reversing the orientation of each component 
in~$\tilde{T}$. Define~$\epsilon(T)$ by deleting~$\tilde{T}$. Now 
let~$D$ be a unitrivalent diagram~$D$ with a distinguished 
subset~$\tilde{\Gamma}$ of its skeleton components. Define~$d(D)$ 
by replacing each skeleton component in~$\tilde{\Gamma}$ by two 
copies, and by summing over all ways of  lifting vertices of~$D$ 
that lie on~$\tilde{\Gamma}$ to the new skeleton. Define~$s(D)$ by 
reversing the orientation of the components in~$\tilde{\Gamma}$ 
and by multiplying with~$\prod_{C\in\tilde{\Gamma}}(-1)^{n_C}$ 
where~$n_C$ is the number of vertices lying on the skeleton 
component~$C$ of~$D$. If $n_C>0$ for some component~$C$ 
of~$\tilde{\Gamma}$, then define~$\epsilon(D)=0$. 
Define~$\epsilon(D)$ by deleting the components 
in~$\tilde{\Gamma}$ in the remaining case. The 
composition~$T_1\circ T_2$ of non-associative tangles~$T_1, T_2$ 
with $\source(T_1)=\target(T_2)$ is defined by placing $T_1$ on 
the top of~$T_2$. For diagrams~$D_i$ in~$\A(\Gamma(T_i))$ a 
composition~$D_1\circ D_2$ is defined similarly. In the following 
lemma we state generalizations of well-known properties of~$Z$. 

\begin{lemma}\label{l:A0gen}
Let $T$, $T_1$, $T_2$ be non-associative tangles with dotted 
circles. 

(1) Assume that some of the components of~$T$ without dots are 
distinguished. Then we have\footnote{As in~\cite{LM3} we must 
assume for the first property of~$\Aarhus_0$ that an even 
associator is used in the definition of~$Z$. This causes no 
restrictions in Theorems~\ref{t:ZlmoNabla} and~\ref{t:nablacan} 
because for links~$L$, the invariants~$\Aarhus(L)$ 
and~$Z^{LMO}(L)$ do not depend on the choice of an associator.}  

$$ d(\Aarhus_0(T))=\Aarhus_0(d(T)) \quad ,\quad 
s(\Aarhus_0(T))=\Aarhus_0(s(T))\quad,\quad 
\epsilon(\Aarhus_0(T))=\Aarhus_0(\epsilon(T)). $$ 

(2) Assume that~$\source(T_1)=\target(T_2)$. Then 

$$ \Aarhus_0(T_1\circ T_2)=\Aarhus_0(T_1)\circ\Aarhus_0(T_2).$$ 

(3) We have

$$ \Aarhus_0(T)=\bar{\chi}_Y((\exp(P)),$$

where $P$ is a series of connected diagrams in~$\A(\emptyset,Y)$ 
and~$Y$ is a set in bijection with the components of~$T$ without 
dots. 
\end{lemma}

The proof of Lemma~\ref{l:A0gen} is straightforward. Statements 
similar to Lemma~\ref{l:A0gen} hold for~$\Aarhus(L)$. 
 
\section{Proofs of Theorems~\ref{t:ZlmoNabla} and~\ref{t:nablacan}}\label{s:proofs}

Recall from Equation~(\ref{e:nablacan}) that for links~$L$ 
in~$S^3$ we have~$c\nabla(L)_{\vert\;t^{1/2}=e^{h/2}}= 
W_\nabla\circ Z(L)$ with $c=h/({e^{h/2}-e^{-h/2}})$. 
Equation~(\ref{e:nablacan}) is proven in~\cite{LM1} and~\cite{BNG} 
by showing that~$W_\nabla\circ Z$ satisfies a skein relation 
and~$W_\nabla\circ Z(O)=c$. With the methods of this proof one can 
show directly that~$W_\nabla\circ\Aarhus$ satisfies the same skein 
relation for links in a rational homology sphere 
and~$W_\nabla\circ\Aarhus(O)=c$, but this does not imply 
Theorem~\ref{t:nablacan}. In this section we present a proof of 
Theorem~\ref{t:nablacan} based on Equation~(\ref{e:nablacan}). 
Then we prove Theorem~\ref{t:ZlmoNabla} by using 
Theorem~\ref{t:nablacan}. 

Let $L=L'\cup L''$ be a framed link in a rational homology 
sphere~$M$. Denote the components of~$L'$ (resp.~$L''$) by~$K_x$ 
with $x\in X'$ (resp.~$x\in X''$) and their framings by~$\mu_x$. 
For $x,y\in X'\cup X''$ let $l_{xy}=\lk(\mu_x,K_y)$ be linking 
numbers in~$M$, let the submatrix corresponding to~$L'$ be 
invertible and denote its inverse by~$(l^{xy})_{x,y\in X'}$. In 
the following lemma we recall how the linking numbers transform 
under surgery. 
 
\begin{lemma}\label{l:linkgen}  For $i,j\in X''$ the linking 
numbers $\tilde{l}_{ij}=\lk(\mu_i,K_j)$ of $L''\subset M_{L'}$ are 
given by 

$$ \tilde{l}_{ij}=l_{ij}-\sum_{x,y\in X'}l_{ix}l^{xy}l_{yj}. $$ 
\end{lemma} \begin{proof} Denote the meridians of the components 
of~$L$ by~$m_x$. In $H_1(M\setminus (L'\cup L''),\Q)$ the framings 
$\mu_y$ can uniquely be expressed as $\mu_y=\sum_{j\in X'\cup X''} 
l_{yj}m_j$. This implies for $x\in X'$ that 

$$\sum_{y\in X'} l^{xy} \mu_y=m_x+\sum_{y\in X', j\in 
X''}l^{xy}l_{yj}m_j.$$ 

In $H_1(M_{L'}\setminus L'',\Q)=H_1(M\setminus (L'\cup 
L''),\Q)/(\mu_x)_{x\in X'}$ we obtain the following unique 
expression of~$\mu_i$ ($i\in X''$) in terms of the meridians~$m_j$ 
($j\in X''$) of~$L''\subset M_{L'}$:

$$\mu_i=\sum_{j\in X'\cup X''} l_{ij}m_j=\sum_{j\in X''} 
l_{ij}m_j-\sum_{x\in X',y\in X',j\in X''}l_{ix}l^{xy}l_{yj}m_j.$$ 
This implies the lemma. \end{proof}

The following lemma tells us that the linking numbers of a 
link~$L\subset M$ can be recovered from~$\Pstrut(\Aarhus(L))$. 

\begin{lemma}\label{l:linkA} Let~$L$ be a link with integral 
framing in a rational homology sphere~$M$. Let the components 
of~$L$ be in bijection with a set~$X$. Let 
$(\tilde{l}_{ij})_{i,j\in X}$ be the linking matrix of~$L$. Then 

$$\Pstrut(\Aarhus(L))=\exp\left( \frac{1}{2}\sum_{i,j\in 
X}\tilde{l}_{ij}\strutn{i}{j}\right).$$ \end{lemma} \begin{proof} 
Choose a diagram of $L'\cup L''\subset S^3$ such that 
$(S^3_{L'},L'')\cong (M,L)$ and put dots on the components of 
$L'$. Let $(l_{xy})_{x,y\in X'\cup X''}$ be the linking matrix of 
$L'\cup L''$ and let $(l^{xy})_{x,y\in X'}$ be the inverse of the 
linking matrix of~$L'$. Then for a series~$P$ (resp.\ $\tilde{P}$) 
of diagrams in~$\A(\emptyset, X'\cup X'')$ (resp.\ 
in~$\A(\emptyset, X'')$) that contains no struts and has 
degree-$0$-term~$1$, we have 
 
\begin{eqnarray*} \cZ(L'\cup L'') & = & 
\chib_{X''}\left(P\disuni\exp\left(\frac{1}{2}\sum_{x,y\in X'\cup 
X'' }l_{xy}\strutn{x}{y}\right)\right)\quad\mbox{and} 
\end{eqnarray*} \begin{eqnarray*} \Aarhus_0(L'\cup L'') & =& 
\chib_{X''}\left(<P\disuni\exp\left(\frac{1}{2}\sum_{i,j\in X'' 
}l_{ij}\strutn{i}{j}+\sum_{i\in X'', x\in X' 
}l_{ix}\strutn{i}{x}\right)\right., \\ 
&&\left.\qquad\qquad\exp\left(-\frac{1}{2}\sum_{x,y\in X' 
}l^{xy}\strutu{\partial x}{\partial y}\right)>\right)\\ &=& 
\chib_{X''}\left(\tilde{P}\disuni 
\exp\left(\frac{1}{2}\sum_{i,j\in 
X''}l_{ij}\strutn{i}{j}-\frac{1}{2}\sum_{i,j\in X'', x,y\in X' 
}l_{ix}l^{xy}l_{yj}\strutn{i}{j}\right)\right). \end{eqnarray*} 

Since~$\Pstrut(\Aarhus(L))=\Pstrut(\Aarhus_0(L'\cup L''))$, 
Lemma~\ref{l:linkA} follows from Lemma~\ref{l:linkgen}. 
\end{proof} 

For technical reasons we fix a representative of each 
homeomorphism-class of connected compact surfaces with boundary. 
We call this representative~$\Sigma$ a standard surface and equip 
it with a decomposition into a single vertex~$v\cong I\times I$ 
(also called coupon) with bands~$B_i\cong I\times I$ that are 
glued along~$I\times\{0,1\}$ to the upper boundary~$I\times 
\{1\}$~of~$v$. Call the part~$I\times \{0\}$ of~$v$ its 
distinguished lower boundary. We orient the cores $I\times 
\{1/2\}$ of the bands~$B_i$ ($i=1,\ldots,\rank\, H_1(\Sigma)$). An 
example is shown on the left side of 
Figure~\ref{f:standardsurface}. 

\begin{figure}[!h]
\centering \setbox1=\hbox{\input{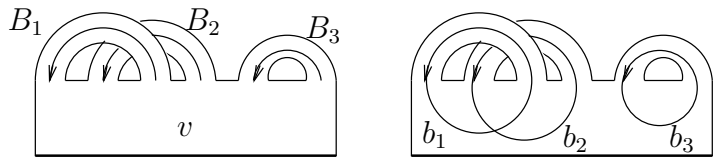}}
$\vcenter{\box1}$ \caption{A standard surface~$\Sigma$ and 
a basis~$(b_i)$ of~$H_1(\Sigma)$}\label{f:standardsurface} 
\end{figure}

We associate a basis of~$H_1(\Sigma)$ to the ribbon graph 
decomposition of~$\Sigma$ as shown in 
Figure~\ref{f:standardsurface} by an example. The orientation 
of~$b_i$ is determined by the orientation of the core of the 
band~$B_i$. An embedding of a standard surface into~$R^2\times I$ 
is an example of a ribbon graph in the sense of Section~8 
of~\cite{KaT}. Ribbon graphs without vertices can canonically be 
identified with framed tangles. We will use this identification in 
the following. 

From now on we use the term Seifert matrix of a Seifert 
surface~$\Sigma\subset M$ always with respect to a basis 
of~$H_1(\Sigma)$ obtained by identifying~$\Sigma$ with a standard 
surface in some freely chosen way. We use the same basis for a 
matrix of the intersection form of~$\Sigma$. 

\begin{keylemma}\label{kl:AKS} Let~$K$ be a knot in a rational 
homology sphere~$M$ bounding a Seifert surface~$\Sigma$. Let~$V$ 
be a Seifert matrix of~$\Sigma$. Then the power 
series~$W_\nabla\circ\Aarhus(K)$ depends only on~$V$. The 
coefficient of~$h^i$ in this series is a polynomial in the entries 
of~$V$. \end{keylemma} 

Let us prepare the proof of Keylemma~\ref{kl:AKS}. We will make 
some statements more generally for links instead of knots. 
Let~$V=(v_{ij})$ be a Seifert matrix. Choose a null-homologous 
link~$L$ in a rational homology sphere~$M$ with Seifert 
surface~$\Sigma$ and Seifert matrix~$V$. The homeomorphism type 
of~$\Sigma$ is determined by the similarity type of~$V-V^*$. There 
exists a link with integral framing~$\tilde{L}\subset M$ such 
that~$M_{\tilde{L}}=S^3$. The link~$\tilde{L}$ can be chosen to be 
disjoint from~$\Sigma$ because changing crossings 
between~$L=\partial\Sigma$ and~$\tilde{L}$ preserves the property 
that~$M_{\tilde{L}}=S^3$. Therefore~$\Sigma\subset M$ can be 
obtained from a surface~$\Sigma''\subset S^3$ by surgery along a 
link~$L'\subset S^3\setminus\Sigma''$. The identification 
of~$\Sigma$ with a standard surface induces an identification 
of~$\Sigma''$ with a standard surface. In a diagram of~$L'\cup 
\Sigma''$ the vertex~$v$ of~$\Sigma''$ can be pulled downwards, 
such that the diagram~$L'\cup\Sigma''$ is of the form~$(L'\cup 
T_1'')\circ T_2''$ where~$T_2''$ is a planar diagram of a 
neighborhood of~$v$ and the distinguished lower boundary of~$v$ is 
the lowest part of the diagram. Put dots on the components 
of~$L'$. An example is shown in Figure~\ref{f:ribbongr}. 

\begin{figure}[!h]
\centering \setbox1=\hbox{\input{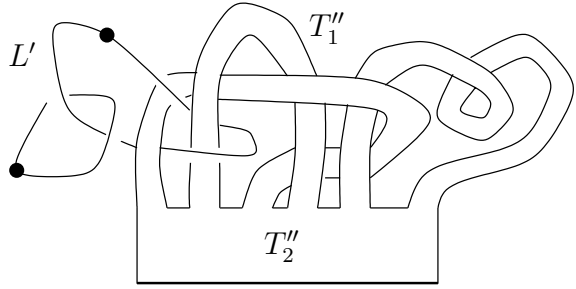}}
$\vcenter{\box1}$ \caption{A diagram of $L'\cup \Sigma''$
}\label{f:ribbongr} 
\end{figure}

The components of~$T_1''$ are in bijection with the set 
$X''=\{1,\ldots, \rank\; H_1(\Sigma)\}$. Regard~$T_1''$ as a 
non-associative framed tangle with parentheses of the 
form~$((((...).).).)$ on~$\source(T_1'')$. Let $F=(f_{ij})=V-V^*$ 
be the matrix of the intersection form of~$\Sigma$ and 
let~$U=(u_{ij})=1/2(V+V^*)=V-1/2F$. 
 
\begin{lemma}\label{l:T1}
With the notation from above we have \nopagebreak{}
 $$ \Pstrut 
(\Aarhus_0 (L'\cup T_1''))=\exp\left(\frac{1}{2}\sum_{i,j\in X''} 
u_{ij}\strutn{i}{j}\right). $$ \end{lemma} \begin{proof} Let 
$K_i^+$ (resp.~$K_i^-$) be a knot in the upper part 
$\Sigma^+\subset M$ (resp.\ in the lower part~$\Sigma^-\subset M$) 
of a tubular neighborhood of~$\Sigma$ representing the $i$-th 
basis element of~$H_1(\Sigma)$. Let the knot $K_i^-\cong K_i^+$ 
have the framing~$\lk(K_i^-,K_i^+)=v_{ii}$ induced by the 
surface~$\Sigma$. First consider~$i\not=j\in X''$. Define~$P_{ij}$ 
as the composition of~$\Pstrut$ with the projection to the part 
containing only powers of the strut~$\strutn{i}{j}$. 
Lemma~\ref{l:linkA} implies that $P_{ij}(\Aarhus(K_i^-\cup 
K_j^+))=\exp({v_{ij}\strutn{i}{j}})$. Represent~$K_i^-\cup 
K_j^+\subset M$ by a surgery diagram~$(L'\cup S_1'')\circ S_2''$ 
where the tangle~$S_1''$ consists of the $i$-th and $j$-th framed 
strands of~$T_1''$ and~$S_2''$ is a~$0$-framed tangle consisting 
of two intervals close to~$T_2''$. See Figure~\ref{f:LSi} for an 
example (compare Figures~\ref{f:standardsurface} 
and~\ref{f:ribbongr}). In this figure the dotted line 
separates~$S_2''$ from $L'\cup S_1''$ and is not a part of the 
diagram. 

\begin{figure}[!h] \centering \setbox1=\hbox{\input{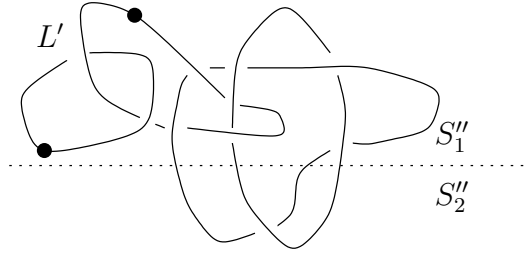}}
$\vcenter{\box1}$ \caption{A surgery 
diagram $(L'\cup S_1'')\circ S_2''$ of $K_1^-\cup K_2^+$ 
}\label{f:LSi} \end{figure}

We have~$\Aarhus_0(S_2'')=Z(S_2'')$ and the explicit description 
of~$Z$ (see~\cite{LM2}) 
implies

$$ \Pstrut(\Aarhus_0(S_2''))=\exp((1/2)f_{ij}\strutn{i}{j}). $$ 
 
Observe the following property of~$P_{ij}$: 

$$P_{ij}\left(\Aarhus_0(L'\cup S_1'')\circ 
\Aarhus_0(S_2'')\right)=P_{ij}(\Aarhus_0(L'\cup S_1''))\disuni 
P_{ij}(\Aarhus_0(S_2'')).$$ 
 
The last two formulas and Part~(2) of Lemma~\ref{l:A0gen} imply 

$$P_{ij}(\Aarhus_0(L'\cup S_1'' 
))=\exp((v_{ij}-f_{ij}/2)\strutn{i}{j})=\exp(u_{ij}\strutn{i}{j}).$$
  
Using Lemma~\ref{l:A0gen} for~$\epsilon$ we see 
that~$P_{ij}(\Aarhus(L'\cup T_1''))=\exp(u_{ij}\strutn{i}{j})$. 
For $i=j$ Lemma~\ref{l:linkA} implies 
$P_{ii}(\Aarhus(K_i^\pm))=\exp((1/2){v_{ii}\strutn{i}{i}})=\exp( 
(1/2) u_{ii}\strutn{i}{i})$.  We apply Lemma~\ref{l:A0gen} 
for~$\epsilon$ as above and obtain~$P_{ii}(\Aarhus_0(L'\cup 
T_1''))=\exp( (1/2)u_{ii}\strutn{i}{i})$. By Part~(3) of 
Lemma~\ref{l:A0gen} we have $\Pstrut(\Aarhus_0(L'\cup 
T_1''))=\Disuni_{i\leq j}P_{ij}(\Aarhus_0(L'\cup T_1''))$ which 
completes the proof.  \end{proof}

Starting from~$V$ we made a lot of choices in the definition 
of~$T_1''$. Since~$\Aarhus(L)$ is an invariant only the choice 
of~$L\subset M$ can influence~$W_\nabla\circ\Aarhus(L)$. Now we 
are ready to show that for knots~$L$ the 
invariant~$W_\nabla\circ\Aarhus(L)$ depends only on~$V$. 

\medskip

\begin{proof}[Proof of Keylemma~\ref{kl:AKS}] We use the notation 
from 
above. 
For suitable distinguished components of~$T_1''$ and of~$d(T_1'')$ 
the tangle~$s(d(T_1''))$ coincides with the part of the framed 
oriented boundary of~$\Sigma''$ that belongs to~$T_1''$. 
Let~$T_3''$ be the part of the framed oriented boundary 
of~$\Sigma''$ that belongs to~$T_2''$. We regard~$T_3''$ as a 
non-associative tangle with $\target(T_3'')=\source(s(d(T_1'')))$. 
The invariant $Z(T_3'')=\Aarhus_0(T_3'')$ depends only on $\rank\; 
H_1(\Sigma)$. Since we know that the Seifert matrix~$V$ is chosen 
with respect to a basis induced by a standard surface, the 
definition of the map~$\chi_{X''}:\A(\emptyset,X'')\longrightarrow 
\A(\Gamma(T_1''),\emptyset)$ depends only on~$V$ (see 
Equation~(\ref{e:defchi})). We will show below that for knots~$L$ 
all terms in~$\Aarhus_0(L'\cup T_1'')$ that contain an internal 
vertex do not contribute to~$W_\nabla(\Aarhus(L))$. 
Equation~(\ref{e:A0A}), Lemma~\ref{l:A0gen} and Lemma~\ref{l:T1} 
then imply 

\begin{eqnarray*} W_\nabla(\Aarhus(L)) & = & 
W_\nabla(\Aarhus_0(L'\cup\partial\Sigma''))\\ & = & 
W_\nabla(s(d(\Aarhus_0(L'\cup T_1'')))\circ Z(T_3''))\\
 & = & W_\nabla\left(s(d(\chi_{X''}(\Pstrut(\Aarhus_0(L'\cup T_1'')))))\circ 
Z(T_3'')\right)\\ & = & W_\nabla\left(s\circ 
d\circ\chi_{X''}\left(\exp\left(\frac{1}{4}\sum_{i,j\in X''} 
(v_{ij}+v_{ji})\strutn{i}{j}\right)\right)\circ Z(T_3'')\right). 
\end{eqnarray*} 

This will show that~$W_\nabla(\Aarhus(L))$ is determined by the 
Seifert matrix~$V$. Obviously, the coefficients of~$h^i$ 
in~$W_\nabla(\Aarhus(L))$ are polynomials of degree~$\leq i$ in 
the entries of~$V$. This will prove the keylemma.  

It remains to consider diagrams~$D$ in~$\Aarhus_0(L'\cup T_1'')$ 
with an internal vertex~$u$. In~$s(d(D))$ each of the edges 
incident to~$u$ is either connected to another internal vertex or 
appears twice, namely as the difference of the two ways of lifting 
it to the skeleton~$\Gamma(s(d(T_1'')))$. We represent this 
difference by a box in Figure~\ref{f:3neighbors}. 

\begin{figure}[!h]
\centering \setbox1=\hbox{\input{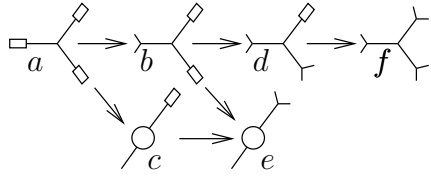}}
$\vcenter{\box1}$ \caption{Replacing differences of 
univalent vertices by internal vertices}\label{f:3neighbors} 
\end{figure}

A neighborhood of the internal vertex~$u$ looks like in one of the 
possibilities (a)-(f) in Figure~\ref{f:3neighbors}. When we push a 
lifted vertex in the box along the 
circle~$\Gamma(\partial\Sigma'')$, then it will finally cancel 
with the second lifted vertex. By the (STU)-relation we can 
replace a box in Figure~\ref{f:3neighbors} by a sum of diagrams 
with an additional internal vertex. More precisely, a part of the 
diagram looking like in (a), (b), (c), or (d) in 
Figure~\ref{f:3neighbors} is replaced by a sum of diagrams where a 
neighborhood of~$u$ looks like in diagrams that can be reached by 
following a directed arrow in Figure~\ref{f:3neighbors}. When we 
apply this procedure to all boxes, we will finally end up with 
possibilities~(e) and~(f). By Lemma~\ref{l:Wnabla} all diagrams 
that have a subdiagram as in (e) or~(f) are sent to~$0$ 
by~$W_\nabla$. \end{proof} 

Let us recall a fact about knots (and links) in $S^3$ (see 
Proposition~8.7 of~\cite{BuZ}). 
 
\begin{fact}\label{f:sm_link}
Let~$V$ be a $n\times n$-matrix over~$\Z$ such that~$V-V^*$ is a 
matrix of the intersection form of a surface. Then~$V$ is a 
Seifert matrix of a link in~$S^3$. 
\end{fact}

Since for all Seifert forms~$s$ the intersection form of~$\Sigma$ 
is equal to~$s-s^*$, we see that Seifert forms of a fixed 
surface~$\Sigma$ are a subset of an affine space whose 
associated~$\Q$-vector space are symmetric $\Z$-bilinear forms 
on~$H_1(\Sigma)$ with values in~$\Q$. By Fact~\ref{f:sm_link} 
Seifert forms of Seifert surfaces~$\Sigma$ in~$S^3$ are a lattice 
of full rank in this affine space. 

\medskip

\begin{proof}[Proof of Theorem~\ref{t:nablacan}] By 
Proposition~\ref{p:nablaKS} and Keylemma~\ref{kl:AKS} the 
coefficients of $h^i$ in the two power series 
$\frac{h}{e^{h/2}-e^{-h/2}}\nabla(K)_{\vert t^{1/2}=e^{h/2}}$ and 
$W_\nabla\circ\Aarhus(K)$ only depend on a Seifert matrix~$V$ of a 
knot~$K$ and are polynomials~$p_i$ and~$q_i$ in the entries 
of~$V$. By Equation~(\ref{e:nablacan}) we have~$p_i(V)=q_i(V)$ for 
all Seifert matrices of knots in~$S^3$. Fact~\ref{f:sm_link} 
implies that~$p_i=q_i$ for all~$i$. 
\end{proof}

The following lemma is a straightforward extension of a result 
of~\cite{GaH}. 

\begin{lemma}\label{l:ZlmoKMK} Let~$K$ be a $0$-framed knot in a 
rational homology sphere~$M$. Then any of the series 
$Z^{LMO}(M_K)$, $W_\nabla\circ Z^{LMO}(K)$, 
$W_\nabla\circ\Aarhus(K)$ can be computed from any other of these 
series. \end{lemma} \begin{proof}[Sketch of proof] The invariants 
$Z^{LMO}$ and $\Aarhus$ of $K\subset M$ differ only by 
normalization (see Equation~(\ref{e:ALMO})). Let $C=Z^{LMO}(K)$. 
Then we have $Z^{LMO}(M_K)=\mbox{$<\sigma(\nu\# C)>$}$ 
with~$\nu=Z(O)$.  The following four steps show that~$W_\nabla(C)$ 
can be calculated from~$Z^{LMO}(M_K)$ and vice versa. This will 
complete the proof. 

1) $W_\nabla(C)$ depends only on the wheel-part~$\Pwh(C)$ of~$C$ 
(Lemma~\ref{l:WnablaPwh}). 

2) $\Pwh(C)$ can be calculated from~$W_\nabla(C)$ 
because~$\Pwh(C)=\exp (P)$ where~$P$ is a formal series of 
connected wheels (see Part~(3) of Lemma~\ref{l:A0gen}), 
$W_\nabla(C)=\exp(W_\nabla(P))$, and~$W_\nabla$ is injective on 
connected wheels (Lemma~\ref{l:WnablaPwh}). 

3) $\sigma(\nu\# C)$ contains no struts because~$K$ is $0$-framed 
(see Lemma~\ref{l:linkA} and Equation~(\ref{e:ALMO})). All 
remaining non-vanishing diagrams in~$\A(\emptyset,\{x\})$ have at 
least as many internal vertices as univalent vertices. This 
implies that~$Z^{LMO}(M_K)$ depends only on~$\Pwh(\nu\# 
C)=\Pwh(\nu)\disuni\Pwh(C)$. 

4) The map $<\cdot>$ is injective on wheels (see \cite{GaH}, 
Lemma~3.1, use the ${\rm sl}_2$-weight system on~$\A(\emptyset)$ 
to see that~$<\cdot>$ is injective on connected wheels). Therefore 
$\Pwh(\nu)\disuni\Pwh(C)$ can be calculated 
from~$Z^{LMO}(M_K)=<\Pwh(\nu)\disuni\Pwh(C)>$. $\Pwh(\nu)$ is 
invertible. \end{proof} 

Now we prove the main result of this paper. 

\medskip

\begin{proof}[Proof of Theorem~\ref{t:ZlmoNabla}]  By 
Lemmas~\ref{l:surga} and~\ref{l:ZlmoKMK} and by 
Equation~(\ref{e:nablaKMK}) it is sufficient to show that for a 
null-homotopic knot~$K$ in a rational homology sphere each of the 
invariants~$\nabla(K)$ and~$W_\nabla\circ\Aarhus(K)$ can be 
computed from the other one. This statement follows from 
Theorem~\ref{t:nablacan}. \end{proof}

} 
\end{document}